\documentclass[12pt]{article}
\newtheorem{ttt}{Theorem}
\newtheorem{ddd}{Definition}[section]

\newtheorem{lll}{Lemma}[section]

\newcommand{\aut}{automorphism\hspace{2pt}}

\title{The group of automorphisms of semigroup of endomorphisms
of free commutative and free associative algebras}
\author{A. Berzins\\University of Latvia\\ e-mail: aberzins@latnet.lv}
\date{August 2004}

\begin{document}
\maketitle

AMS \textit{Mathematics Subject Classification}: 14A99, 14P05


\begin{abstract}

In this paper are described the groups of automorphisms of
semigroup End$(W(X))$, where $W(X)$ is free commutative or free
associative algebra.

\end{abstract}

\section{Introduction}

The basis of the classical algebraic geometry is the Galois
correspondence between \textit{P}-closed ideals in
$P[x_1,x_2,\ldots,x_n]$ and algebraic sets in the affine space
$P^n$. Having noted that a point $a=(a_1,a_2,\ldots,a_n)$ unique
determines a homomorphism $x\in\textrm{Hom}(P[X])$
($s(x_i)=a_i$), B. Plotkin in [11-13] makes foundation of the
universal algebraic geometry, i.e., algebraic geometry for
arbitrary variety $\Theta$ of universal algebras. Let us fix a
variety $\Theta$ and an algebra $H\in\Theta$. Let $W=W(X)$ be a
free algebra in $\Theta$ with finite set of generators. The set
$\textrm{Hom}(W,H)$ we consider as affine space of points over
$H$. Then arbitrary congruence $T$ in $W$ determines the set of
points in affine space $\textrm{Hom}(W,H)$:
$$T'_H=A=\{\mu\in\textrm{Hom}|T\subset\textrm{Ker}\mu\}.$$
A set of points $A\subset\textrm{Hom}(W,H)$ determines congruence
of W
$$A'=A'_W=T=\bigcap_{\mu\in A}\textrm{Ker}\mu.$$
We call the set of points $A$ such that $A=T'$ for some $T$ an
algebraic set in $\textrm{Hom}(W,H)$. A relation $T$ with $T=A'$
for some $T$ is a congruence in $W$. We call such a congruence an
$H$-closed one.

It was found [3, 4, 11-14] that many problems in the universal
algebraic geometry, such as geometric equivalence, geometric
similarity, isomorphism and equivalence of categories of algebraic
sets and varieties depend on structure of ${\rm Aut}(\Theta^0)$
and Aut(End($W$)), where $\Theta^0$ is the category of algebras
$W=W(X)$ with the finite $X$ that are free in $\Theta$, and $W$ is
a free algebra in $\Theta$. Structure of ${\rm Aut(Com}-P))^0$
(i.e, the classical case) was described in [2]. Here was
introduced variant of the concept of quasiinner automorphism -
basic  concept for description of Aut(End($W$)) and
$\textrm{Aut}(\Theta^0)$. Later groups $\textrm{Aut}(\Theta^0)$
were described for categories of free Lie algebras and free
associative algebras [1, 5, 8, 9]and for categories of free
modules and free Lie modules [5, 6]. Some other definition of
quasiinner automorphism was used in mentioned papers. Also are
described ${\rm Aut}(\Theta^0)$ for varieties of groups,
semigroups and some another varieties. For some varieties is
described Aut(End($W$)), but these problems are more difficult.
In particular, this problem was not solved even in the classical
case of free commutative algebras, i.e. the structure of
Aut(End$P[x_1,\ldots ,x_n]$)) is not describe if $n>2$. In all
solved cases groups Aut(End($W$)) are generated by semi-inner
automorphisms and mirror automorphism. In [3] this problem in
classical case was solved, if the field $P$ is algebraically
closed.

I want to note that the proof in [3] does not use the description
of the group ${\rm Aut}(P[X])$). It is very important, because
structure of the group ${\rm Aut}(P[x_1,\ldots ,x_n])$ is not
described, when $n>2$, but counterexamples show that structures of
${\rm Aut}(P[x_1,\ldots ,x_n])$ and ${\rm Aut}(P[x,y]$) are
principally different. Now in the presented paper the group
Aut(End($W_n$)) is described for commutative and associative
algebras over arbitrary infinite field.

\section{Definitions}

Let us recall some definitions for variety of commutative algebras
$\textrm{Com}-P$ (see \textrm{[4,13]} for general case).

\begin{ddd}

Let $P[X]=P[x_1,\ldots ,x_n]$ be a free commutative algebra over
a field $P$ with a finite set of generators and
$\tau\in\textrm{Aut}(\textrm{End}(P[X]))$. It is known
$\textrm{[4]}$ that there exists a bijection $\mu:P[X]\mapsto
P[X]$ such that for every $s\in \textrm{End}(P[X])$
$$ s^{\tau}=\mu s\mu^{-1}.$$
Such representation of $\tau$ is called a representation of $\tau$
as quasiinner automorphism generated by $\mu$.

\end{ddd}
Of course, arbitrary set of substitutions does not generate an
automorphism, But there are two kinds of bijections $\mu$ in
commutative case and three in the associative case that
generate an automorphism.\\

1. $\mu=\overline{\alpha}$, where $\overline{\alpha}$ is the
natural extension of an automorphism $\alpha \in
\textrm{Aut}(P[X])$ on $W$, i.e. $\overline{\alpha}\left( \sum
a_ku_k\right)=\sum \alpha(a_k)u_k$.
 We shall write $\overline{\alpha}=\alpha$ and call it "\aut of
 field".\\

 2. $\mu=\eta\in \textrm{Aut}(P[X])$.\\

\begin{ddd}

The automorphism $\tau$ generated by an automorphism
$\mu\in\textrm{Aut}(P[X])$ is called inner automorphism.

\end{ddd}

\begin{ddd}

Note, that every $\alpha \in \textrm{Aut}(P)$ belongs to
normalizer of subgroup  $\textrm{Aut}(P[X]))$ in the group of all
bijections of $P[X]$. So, every product of elements of
$\textrm{Aut}(P)$ and $\textrm{Aut}(P[X])$ my be represented in
the form $\mu=\alpha\eta$, $\alpha \in \textrm{Aut}P$,
$\eta\in\textrm{Aut}(P[X])$. Automorphism
$\tau\in\textrm{Aut}(\textrm{End}(P[X])$ generated by $\mu$ is
called semiinner.

\end{ddd}

3. Now we describe the third type of $\mu$ for associative
algebras. Let $S=S(X)$ be a free semi-group. For every
$u=x_{i_1}x_{i_2}\ldots x_{i_n}$ in $S$, take
$\overline{u}=x_{i_n}\ldots x_{i_2}x_{i_1}$. Then $u\rightarrow
\overline{u}$ is an anti\aut of the semigroup S.

\begin{ddd}

Let $W=W(X)$ be a free associative algebra. For every its element
$\omega=\lambda_0+\lambda_1u_1+\cdots+\lambda_ku_k$ denote
$\beta(\omega)=\overline{\omega}=\lambda_0+\lambda_1\overline{u_1}+\cdots+
\lambda_k\overline{u_k}$. The transition $\beta:W(X)\rightarrow
W(X)$ is an antiautomorphism of the algebra $W$, and it generates
the automorphism $\delta\in \textrm{Aut}(\textrm{End}(W(X))$,
$s^\delta =\beta s \beta^{-1}$. We call $\beta$ mirror
antiautomorphism of $W$ and $\delta$ mirror automorphism of
$\textrm{Aut}(\Theta^0)$.

\end{ddd}

Here $\delta$ is not inner and is not semi-inner, but is
quasi-inner. Note that

\begin{enumerate}

\item $\delta$ belongs to the normalizer of subgroup
$\textrm{Sinn}(\textrm{End}(W(X))$ in
$\textrm{Aut}(\textrm{End}(W(X))$ that consists of all semi-inner
automorphisms;

\item every antiautomorphism of $W$ is a product of the mirror
antiautomorphism $\beta$ and an automorphism of $W$.

\end{enumerate}

Let now $W(X)$ be a free finitely generated algebra in arbitrary
variety $\Theta$.

\begin{ddd}

Bijection $\mu : W(X)\mapsto W(X)$ is called central, if for every
$s\in \textrm{End}(W(X))$
$$\mu s=s\mu.$$
Algebra $W(X)$ is called central, if every its central bijection
is identical.

\end{ddd}

\begin{ttt}

Let $W(X)$ be a free central algebra and $\mu$ some bijection of
$W(X)$ generating an automorphism
$\tau\in\textrm{Aut}(\textrm{End}(W(X)))$. Then $\mu$ transforms
every base of algebra $W(X)$ to a base of algebra.

\end{ttt}
\textbf{Proof.} [see 4]\\[6pt]

For free algebras in a variety $\Theta$ of algebras over field $P$
(Com-P, Ass-P, and others) we shall correct the definition of
central algebra. Clearly, in every such algebra $W(X)$ we have a
bijection $\mu(u)=au+b$, where $a,b\in P,$ $a\neq 0$ (linear
bijection). This bijection commutate with every endomorphism of
algebra $W(X)$ and, obviously, transforms every base of $W(X)$ to
a base.

\begin{ddd}

Finitely generated algebra $W(X)$ over field $P$ is called almost
central, if every its central bijection is linear.

\end{ddd}

\begin{ttt}

Free finitely generated commutative, associative are almost
central.

\end{ttt}
\textbf{Proof.} Let $\mu$ be a central bijection of $W(X)$,
$\mu(x_i)=r_i(x_1,x_2,\ldots,x_n)\in W(X)$ , and $s\in
\textrm{End}(W(X)$, $s(x_1)=x_1$, $s(x_i)=0$ for $i=1$. Then
$$\mu s(x_1)=\mu (x_1)=r_1(x_1,x_2,\ldots,x_n),$$
$$s\mu(x_1)=s(r_1(x_1,x_2,\ldots,x_n))=r_1(x_1,0,\ldots,0)=r(x_1).$$
So $r_1(x_1,x_2,\ldots,x_n)=r(x_1)$ is a polynomial of one
variable.

For arbitrary $u\in W(X)$ take $s\in \textrm{End}(W(X))$
$s(x_1)=u$, $s(x_i)=x_i$ for $i>1$. Then
$$\mu(u)=\mu s(x_1)=s\mu(x_1)=s(r(x_1))=r(u),$$
and so $\mu(u)=r(u)=a_0+a_1u+\cdots+a_ku^k$. Because $\mu$ is
bijection, so $r$ is linear, Q.E.D..\\[3pt]
\textbf{Corollary.} \textit{From the theorems 1 and 2 we have that
in mentioned algebras every bijection $\mu$, which generate an
automorphism of $\textrm{End}(W(X))$  transform every base of
algebra $W(X)$ to a base of this algebra.}

\section{Linearity of automorphism generating bijection}

Let $W=W(X)$ be a free commutative, associative or Lee algebra
over the field $P$ with a finite set of generators, and $\tau$ be
an automorphism of $\textrm{End}(W(X))$ generated by bijection
$\mu$. We call $s\in \textrm{End}(W(X))$ a constant endomorphism
if $\textrm{Im}(s)=P$. Clearly, if $s$ is constant then $s^{\tau}$
also is constant. Denote the set of all constant endomorphisms by
$Const$.

Let $a\in P$ and $s$ be a constant endomorphism such that
$s(x_1)=a$. Then $\mu(a)=\mu(s(x_1))=s^{\tau}(\mu(x_1))\in P$. So
a restriction of $\mu$ on $P$ is a bijection: $\mu:P\mapsto P$.

Let $\mu(0)=a$ and $\mu(1)=b$. We shall consider the linear
bijection $l:W(X)\mapsto W(X)$, $l(u)=cu+d$ such that $l(a)=0$ and
$l(b)=1$. Since linear bijection generates identical automorphism
of the semigroup $\textrm{End}(W(X))$, so bijections $\mu$ and
$\mu'=l\mu$ generate equal automorphisms. Note that $\mu'(0)=0$
and $\mu'(1)=1$. We denote $\mu=\mu'$ and assume that $\mu(0)=0$
and $\mu(1)=1$.

Since $\mu$ transforms base of algebra $W(X)$ to base, so the
endomorphism $\eta:$ $\eta(x_i)=\mu(x_i)$ is an automorphism of
$W(X)$. Now we shall consider the bijection $\mu'=\eta^{-1}\mu$.
We have for $\mu'$ that $\mu'(x_i)=x_i$.

Now we shall consider the Galois correspondence between ideals in
$W(X)$ and algebraic sets in affine space over $W(X)$. We recall
that every ideal $T$ in $W(X)$ defines an algebraic set of points
in affine space over $W(X)$.
$$T'=T'_{W(X)}=A=\{\eta\in\textrm{End}(W(X))|T\subset\textrm{Ker}\eta\}.$$
A set of points $A\subset\textrm{End}(W(X))$ determines ideal in
$W(X)$:
$$A'=A'_W=\bigcap_{\eta\in A}\textrm{Ker}\hspace{4pt}\eta.$$

An element $u\in W(X)$ is called basic, if it belongs to some base
of algebra $W(X)$.

%lemma 3.1
\begin{lll} If $u$ is an basic element of W(X), then
$$\left\langle u\right\rangle'' =\left\langle u\right\rangle$$
\end{lll}
\textbf{Proof.} Clearly,
$$A=\left\langle u_1\right\rangle'=\{\eta\in\textrm{End}W|\hspace{4pt}\eta(u_1)=0\}$$
and
$$A'=\left\langle u_1\right\rangle''=\bigcap_{\eta\in
A}\textrm{Ker}\hspace{4pt}\eta=\bigcap_{\eta(u_1)=0}\textrm{Ker}\hspace{4pt}\eta.$$
Let $f(u_1,u_2,\ldots,u_n)\in\left\langle u_1\right\rangle''$ and
$\eta:\hspace{4pt}\eta(u_1)=0,\hspace{4pt}\eta(u_i)=u_i$ for
$i>1$. Then $0=\eta(f(u_1,u_2,\ldots,u_n))=f(0,u_2,\ldots,u_n)$,
so $f\in\left\langle u_1\right\rangle$, Q.E.D..

%lemma 3.2
\begin{lll}
Let $\mu$ be an automorphism generating bijection of $W$ such that
$\mu(0)=0,\hspace{2pt}\mu(1)=1$ and $\mu(x_i)=x_i$. Then
$$\left\langle\mu(\alpha
x_1)\right\rangle=\left\langle\mu(x_1)\right\rangle=\left\langle
x_1\right\rangle$$
\end{lll}
\textbf{Proof.} At first note that for any endomorphism
$s\in\textrm{End}(W)$
$$s(x_1)=0\Leftrightarrow s^{\tau}(x_1)=0.$$
Really, if $s(x_1)=0$, then $s^{\tau}(x_1)=s^{\tau}\mu(x_1)=\mu
s(x_1)=0$; and if $s^{\tau}(x_1)=0$, then
$0=s^{\tau}(x_1)=s^{\tau}\mu(x_1)=\mu s(x_1)=0$, i.e. $s(x_1)=0$.

Let $\mu(\alpha x_1)=r_1(x_1,\ldots,x_n)$. Take
$s\in\textrm{End}(W)$ such that $s^{\tau}(x_1)=0$ and
$s^{\tau}(x_i)=x_i$ for $i>1$. Then $0=\mu(0)=\mu s(\alpha
x_1)=s^{\tau}\mu(\alpha
x_1)=s^{\tau}r_1(x_1,\ldots,x_n)=r_1(0,x_1,\ldots,x_n)$, i.e.
$\mu(\alpha x_1)\in\left\langle x_1\right\rangle''=\left\langle
x_1\right\rangle$.

Now let us prove, that $x_1\in\left\langle\mu(\alpha
x_1)\right\rangle$. Note that
$$\langle\mu(\alpha x_1)\rangle=\langle\mu(\alpha
x_1)\rangle''=\bigcap_{\eta^{\tau}:\hspace{2pt}\eta^{\tau}(\mu(\alpha
x_1))=0}\textrm{Ker}\hspace{2pt}\eta^{\tau}.$$ So we must prove
that, if $\eta^{\tau}(\mu(\alpha x_1))=0$ then
$\eta^{\tau}(x_1)=0$. Really,
$$0=\eta^{\tau}(\mu(\alpha x_1))=\mu\eta(\alpha
x_1)=\mu(\alpha\eta(x_1))\Rightarrow\alpha\eta(x_1)=0\Rightarrow\eta^{\tau}(x_1)=0.$$
Q.E.D.

\begin{lll}
Let $\mu$ be an automorphism generating bijection of $W$ such that
$\mu(0)=0,\hspace{2pt}\mu(1)=1$ and $\mu(x_i)=x_i$.Then for every
$a\in P$ and $u\in P[X]$
$$\mu(au)=\mu(a)\mu(u).$$
\end{lll}
\textbf{Proof.} Since $\langle\mu(ax_1)\rangle=\langle
x_1\rangle$, so $\mu(ax_1)=bx_1=b\mu(x_1)$, where $w\in P[X]$.
Substitution $s(x_1)=1$ gives that $b=\mu(a)$, i.e.
$\mu(ax_1)=\mu(a)\mu(x_1)$.

In associative case we consider only the case, where the field
$P$ is of characteristic 0. I must note, that $\mu(ax_1)$ together
with $x_2,\ldots.x_n$ is a base of algebra, and so
$\mu(ax_1)=bx_1+f(x_2,\ldots,x_n)$ (see [7]). Because
$\mu(ax_1)\in\langle x_1\rangle$, so $\mu(ax_1)=bx_1$.

\begin{lll} With the preceding notation
$$\mu(u+v)=\mu(u)+\mu(v).$$
\end{lll}
\textbf{Proof.} Let
$\mu(x_1+x_2)-\mu(x_1)-\mu(x_2)=f(x_1,\ldots,x_n)=g(\mu(x_1),\ldots,\mu(x_n))$,
then $\mu(u_1+u_2)-\mu(u_1)-\mu(u_2)=g(\mu(u_1),\ldots,\mu(u_n)$.
Substitution $s:\hspace{4pt}s(u_1)=u_1,s(u_2)=u_2,s(u_i)=0$ for
$i>2$ gives identity
$$\mu(u_1+u_2)-\mu(u_1)-\mu(u_n)=\overline{g}(\mu(u_1),\mu(u_2)).$$
Clearly, $\overline{g}$ is homogenous of degree 1 of $u_1$ and
$u_2$. Therefore $\overline{g}(u_1,u_2)=au_1+bu_2$ and
$$\mu(u+v)=\mu(u)+\mu(v)+a\mu(u)+b\mu(v).$$
Substitutions $u=0$ and $v=0$ gives that $a=b=0$, Q.E.D..\\[6pt]
\textbf{Corollary.} Note that restriction $\mu$ on $P$ is an
automorphism of $P$ and denote it $\eta$. Let
$\mu'=\eta^{-1}\mu$, then $\mu'$ is identical on $P$. So for
$\mu=\mu'$ we have identity
\begin{equation}
\mu(au)=a\mu(u),
\end{equation}
for every $a\in P$ and $u\in W(X)$.

\section{Description of  $\textrm{Aut}(\textrm{End}(W(X)))$}

Let consider at first the commutative  case.

\begin{lll}
Let $\mu$ be a bijection of $P[X]$ generating an automorphism
$\tau\in \textrm{Aut}(\textrm{End}(P[X]))$ such that $\mu$ is
identical on $P$ and $\mu(x_i)=x_i$. Then for every $u,v\in P[X]$
it holds that
$$\mu(u\cdot v)=\mu(u)\cdot\mu(v)$$
\end{lll}
\textbf{Proof.} Let $s$ be a constant endomorphism of $P[X])$.Then
$s^{\tau}$ also is constant and so
$$s^{\tau}(\mu(x_1x_2)-\mu(x_1)\mu(x_2))=s^{\tau}\mu(x_1x_2)-s^{\tau}\mu(x_1)\cdot
s^{\tau}\mu(x_2)=$$
$$\mu s(x_1x_2)-\mu s(x_1)\cdot \mu
s(x_2)=\mu(a_1a_2)-\mu(a_1)\cdot\mu(a_2)=0,$$
$$\mu(x_1x_2)-\mu(x_1)\mu(x_2)\in\bigcap_{s^{\tau}\in
Const}\textrm{Ker}\hspace{2pt}s^{\tau}=0,$$ Q.E.D..

Now we can formulate the main result for commutative case.

\begin{ttt}
Every automorphism of $\textrm{End}(P[X])$ is semiinner.
\end{ttt}
\textbf{Proof.} It follows from lemma 3.4, lemma 4.1 and identity
(1).\\[6pt]

Now let consider associative case.

\begin{lll}
Let $\mu$ be a bijection of $W(X)$ generating an automorphism
$\tau\in \textrm{Aut}(\textrm{End}(W(X)))$ such that $\mu'$ is
identical on $P$ and $\mu(x_i)=x_i$. Then for every $u,v\in W(X)$
it holds that
$$\mu(u\cdot v)=\mu(u)\cdot\mu(v),$$
or for every $u,v\in W(X)$ it holds that
$$\mu(u\cdot v)=\mu(v)\cdot\mu(u),$$
\end{lll}
\textbf{Proof.} Let
$\mu(x_1x_2)-\mu(x_1)\mu(x_1)=f(x_1,\ldots,x_n)=g(\mu(x_1),\ldots,\mu(x_n))=\overline{g}(\mu(x_1),\mu(x_2)$.
Note that $\overline{g}$ is 1-homogenous of $x_1$ and 1-homogenous
of $x_2$. Following the proof of lemma 3.4 we get identity
$$\mu(uv)=a\mu(u)\mu(v)+b\mu(v)\mu(u).$$
Substitution $v=1$ gives equality $a+b=1$. We shall denote
$\mu(u)=\overline{u}$ and calculate $\mu(xxy)$ .
$$\mu(xxy)=a\overline{xx}\cdot\overline{y}+b\overline{y}\cdot\overline{xx}=$$
$$a(a+b)\overline{x}\cdot\overline{x}\cdot\overline{y}+b(a+b)\overline{y}\cdot\overline{x}\cdot\overline{x}=$$
$$a\overline{x}\cdot\overline{x}\cdot\overline{y}+b\overline{y}\cdot\overline{x}\cdot\overline{x},$$
$$\mu(xxy)=a\overline{x}\cdot\overline{xy}+b\overline{xy}\cdot\overline{x}=$$
$$a^2\overline{x}\cdot\overline{x}\cdot\overline{y}+2ab\overline{x}\cdot\overline{y}\cdot\overline{x}+b^2\overline{y}\cdot\overline{x}\cdot\overline{xy}.$$
So $a=a^2,\hspace{4pt}b=b^2,\hspace{4pt}ab=0,\hspace{4pt}a+b=1$
and we have two solutions $a=1,b=0$ or $a=0,b=1$, Q.E.D..

Now we can formulate the main result for associative algebras.

\begin{ttt}
The group $\textrm{Aut}(\textrm{End}(W(X)))$ is generated by
semiinner automorphisms and mirror automorphism.
\end{ttt}
\textbf{Proof.} It follows from lemma 3.4, lemma 4.2 and identity
(1).\\[6pt]

\noindent{\Large\textbf{Acknowledgments}}\\[12pt]
The author is happy to thank professor B. Plotkin for stimulating
discussions of the results.\\[30pt]

\noindent{\Large\textbf{References}}
\begin{enumerate}

\item A. Berzins, {\it The automorphisms of \,{\rm End}{}$K[x]$}, Proc.
Latvian Acad. Sci., Section B, 2003, vol. 57, no. 3/4, pp. 78-81

\item A. Berzins, {\it Geometric equivalence of algebras}, Int. J.
Alg. Comput., 2001, vol. 11, no. 4, pp. 447-456.

\item A. Berzins, {\it The group of automorphisms of semigroup
End($P[X]$},(to appear)

\item A. Berzins, B. Plotkin and E. Plotkin, {\it Algebraic
geometry in Varieties with the Given Algebra of constants}, J.
Math. Sci., New York, 2000, vol. 102, no. 3, pp. 4039-4070

 \item R. Lipyanski, B. Plotkin {\it Automorphisms of categories
 of free modules and free Lie algebras}, ArXiv math.RA/0502212 (10
 Feb 2005) pp. 14.

 \item Y. Katsov, R. Lipyanski, B. Plotkin, {\it Automorphisms of categories
 of free modules, free semimodules  and free Lie modules}, (to appear)

\item L. Makar-Limanov  {\it Algebraically closed skew field}, Journal of Algebra, Vol.
93, 1985, pp. 117-135.

\item G. Mashevitzky, {\it Automorphisms
of categories of free associative algebras.}, Preprint.

\item G. Mashevitzky, B. Plotkin, E. Plotkin, {\it Automorphisms
of categories of free Lie algebras.}, Journal of Algebra, Vol.
282, 2004, pp. 490-512.

\item B. Plotkin, {\it Algebraic logic, varieties of algebras and
algebraic varieties}, Proc. Int. Alg. Conf., St. Petersburg, 1995,
Walter de Gruyter, New York, London, 1996.

\item B. Plotkin, {\it Varieties of algebras and algebraic
varieties}, Israel J. of Mathematics, 1996, vol. 96, no. 2, pp.
511-522.

\item B. Plotkin, {\it Varieties of algebras and algebraic
varieties. Categories of algebraic varieties}, Sib. Adv. Math.,
1997, vol. 7, no. 2, pp. 64-97.

\item B. Plotkin, {\it Algebras with the same (algebraic)
geometry}, Proc. Steklov Inst. Math, Vol. 242, 2003, pp. 165-196.

\item B. Plotkin, G Zhitomirski, {\it On automorphisms of
categories of universal algebras}, ArXiv math.CT/0411408 (18 Nov
2004). pp. 24.

\end{enumerate}

\end{document}